\documentclass{amsart}
\usepackage{amssymb}
\usepackage{amscd}
\usepackage{hyperref}
\newcommand{\caa}{\mathcal{A}}\newcommand{\cac}{\mathcal{C}}\newcommand{\cau}{\mathcal{U}}\DeclareMathOperator{\rk}{rank}
\newcommand{\ov}{\overline}\newcommand{\frF}{\mathfrak{F}}\newcommand{\cax}{\mathcal{X}}\newcommand{\bbq}{\mathbb{Q}}
\DeclareMathOperator{\et}{\text{\'et}}\newcommand{\cao}{\mathcal{O}}\newcommand{\frp}{\mathfrak{p}}
\DeclareMathOperator{\Frob}{Frob}\DeclareMathOperator{\pic}{Pic}\DeclareMathOperator{\ns}{NS}
\DeclareMathOperator*{\res}{Res}\newcommand{\cag}{\mathcal{G}}\newcommand{\frO}{\mathfrak{O}}
\DeclareMathOperator{\gal}{Gal}\DeclareMathOperator{\aut}{Aut}\newcommand{\cay}{\mathcal{Y}}\newcommand{\bbf}{\mathbb{F}}
\DeclareMathOperator{\tr}{Tr}\newcommand{\frA}{\mathfrak{A}}\newcommand{\frB}{\mathfrak{B}}
\DeclareMathOperator{\sm}{sm}\newcommand{\caf}{\mathcal{F}}\DeclareMathOperator{\gl}{GL}\newcommand{\caz}{\mathcal{Z}}
\newcommand{\call}{\mathcal{L}}\newcommand{\bbz}{\mathbb{Z}}\newcommand{\can}{\mathcal{N}}
\DeclareMathOperator{\spec}{Spec}\DeclareMathOperator{\wt}{wt}\DeclareMathOperator{\codim}{codim}
\newcommand{\cah}{\mathcal{H}}\newcommand{\bbg}{\mathbb{G}}\newcommand{\cav}{\mathcal{V}}\DeclareMathOperator{\supp}{supp}
\DeclareMathOperator{\ur}{ur}\DeclareMathOperator{\tor}{tor}\newcommand{\car}{\mathcal{R}}
\DeclareMathOperator{\gsp}{GSp}\DeclareMathOperator{\mt}{MT}
\newcommand{\bbc}{\mathbb{C}}\newcommand{\bbr}{\mathbb{R}}\DeclareMathOperator{\End}{End}
\DeclareMathOperator{\Aut}{Aut} \newcommand{\und}{\underline}
\newtheorem{theorem}{Theorem}\numberwithin{theorem}{section}
\theoremstyle{plain}
\newtheorem{conjecture}[theorem]{Conjecture}\newtheorem{corollary}[theorem]{Corollary}
\newtheorem{lemma}[theorem]{Lemma}\newtheorem{proposition}[theorem]{Proposition}
\newtheorem*{hypothesis}{Hypothesis}
\theoremstyle{definition}\newtheorem{definition}[theorem]{Definition}
\newtheorem{remark}[theorem]{Remark}\numberwithin{equation}{section}
\begin{document}
\title[Rank abelian varieties]{On the rank of abelian varieties \\over function fields}
\author{Am\'{\i}lcar Pacheco}
\address{Universidade  Federal do Rio de Janeiro (Universidade do Brasil), Departamento de  Ma\-te\-m\'a\-ti\-ca Pura, 
Rua Guaiaquil  83, Cachambi, 20785-050 Rio de Janeiro, RJ, Brasil}
\email{amilcar@acd.ufrj.br}
\thanks{This work was partially supported by CNPq research grant 304424/2003-0, Pronex \newline 41.96.0830.00 and 
CNPq Edital Universal 470099/2003-8. I would like to thank Douglas Ulmer for comments on how to treat the case of 
arbitrary ramification, but the conductor prime to the ramification locus, in the case of elliptic fibrations. I would
also like to thank Marc Hindry for comments on the inequality comparing the conductors of $A$ and $A'$. Finally, I also thank the 
referee for his comments and criticisms.}
\date{August 29, 2005}
\begin{abstract}
Let $\cac$ be a smooth projective curve defined over a number field $k$, $A/k(\cac)$ an abelian variety and $(\tau,B)$ 
the $k(\cac)/k$-trace of $A$. We estimate how the rank of
$A(k(\cac))/\tau B(k)$ varies when we take a finite geometrically abelian cover $\pi:\cac'\to\cac$ defined over $k$.
\end{abstract}
\maketitle
\section{Introduction}

Let $\cac$ be a smooth projective irreducible curve defined over a number field $k$ with function field $K:=k(\cac)$. Let $A/K$ be a
non isotrivial  
abelian variety of dimension $d$ and $\phi:\caa\to\cac$ a proper flat morphism defined over $k$ from a smooth projective irreducible 
variety $\caa$ also defined over $k$ whose generic fiber is $A/K$.
Thus, all smooth fibers of $\phi$ will also be abelian varieties. Let $(\tau,B)$ be the 
$K/k$-trace of $A$. A theorem due to Lang and N\'eron 
states that $A(K)/\tau B(k)$ and $A(\ov{k}(\cac))/\tau B(\ov{k})$ are finitely generated abelian groups. 
Let $g$ be the genus of $\cac$, $\frF_A$ the 
conductor divisor of $A$ on $\cac$ and $f_A$ its degree.

Ogg gave a upper bound for the rank of the latter group in the geometric case (cf. \cite[VI, p.19]{ogg67} and also 
\cite[Theorem 2]{ogg62}), i.e., 
\begin{equation}\label{ogg}
\rk\left(\frac{A(\ov{k}(\cac))}{\tau B(\ov{k})}\right)\le 2d(2g-2)+f_A+4\dim(B).
\end{equation}
In particular, this is also an upper bound for the rank of $A(K)/\tau B(k)$. 

In the case where $A=E$ is an elliptic curve, (\ref{ogg}) reduces to
\begin{equation}\label{shi1}
\rk(E(\ov{k}(C)))\le4g-4+f_E,
\end{equation}
which is a result due to Shioda \cite[Corollary 2]{shi92}. 

We consider a finite cover $\pi:\cac'\to\cac$ defined over $k$ which is geometrically abelian with geometric automorphism group 
$\cag:=\aut_{\ov{k}}(\cac'/\cac)$. Let 
$g'$ be the genus of $\cac'$, $K':=k(\cac')$ and $A':=A\times_KK'$. Let $f_{A'}$ be the degree of the conductor $\frF_{A'}$ of $A'$. 
Note that since the extension $K'/K$ is geometric,
the $K'/k$-trace of $A'$ equals the $K/k$-trace of $A$. By (\ref{ogg}) we have the following upper bound for the rank of $A(K')/\tau 
B(k)$ : 
\begin{equation}\label{ogg2}\begin{aligned}
\rk\left(\frac{A(K')}{\tau B(k)}\right)&\le\rk\left(\frac{A(\ov{k}(\cac'))}{\tau B(\ov{k})}\right)\\
&=\rk\left(\frac{A'(\ov{k}(\cac))}
{\tau B(\ov{k})}\right)\le2d(2g'-2)+f_{A'}+4\dim(B).
\end{aligned}\end{equation}

Our goal in this paper is to give an upper bound for the rank of $A(K')/\tau B(k)$. This
upper bound improves (\ref{ogg2}) in the case where $A$ is a Jacobian variety and $\pi$ is unramified. The latter upper bound
 is an extension of Silverman's
result \cite[Theorem 13]{sil02} obtained in the context of elliptic curves. 
In order to state our main result, as well as his, we need to recall Tate's conjecture for divisors on smooth projective algebraic 
varieties $\cax$ defined over number fields $k$.

Let $l$ be a perfect field and $\ov{l}$ an algebraic closure of $l$. For a smooth projective irreducible algebraic variety
$\cax$ defined over $l$ we denote by $H^i(\cax)$ its $i$-th \'etale cohomology group
$H^i_{\et}(\cax\times_l\ov{l},\bbq_{\ell})$.

Let $\cao_k$ be the ring of integers of $k$ and $\frp$ a prime ideal of $\cao_k$. Let $\Frob_{\frp}\in G_k$ be a
Frobenius automorphism corresponding to $\frp$ and $I_{\frp}\subset G_k$ the inertia group of $\frp$ (well defined up
to conjugation). Define the $L$-function
$$
L_2(\cax/k,s)=\prod_{\frp}\det(1-\Frob_{\frp}q_{\frp}^{-s}\,|\,H^2(\cax)^{I_{\frp}})^{-1}.
$$
Let $\pic(\cax)$ be the Picard group of $\cax$, $\pic^0(\cax)$ the subgroup of divisors algebraically equivalent to
zero, $\ns(\cax)=\pic(\cax)/\pic^0(\cax)$ the N\'eron-Severi group of $\cax$ and $\ns(\cax/k)$ the subgroup of divisor
classes of $\ns(\cax)$ which are defined over $k$. The group $\ns(\cax)$ is finitely generated, hence the same holds
for $\ns(\cax/k)$.

\begin{conjecture}(Tate's conjecture, \cite[Conjecture 2]{tate65})\label{tateconj}
$L_2(\cax/k,s)$ has a pole at $s=2$ of order $\rk(\ns(\cax/k))$.
\end{conjecture}

\begin{remark}
\begin{enumerate}
\item This is a special case of Tate's conjecture which concerns algebraic varieties and algebraic cycles.

\item We do not need to use the hypothesis of the existence of a meromorphic continuation of $L_2(\cax/k,s)$ to
$\Re(s)=2$ interpreting the sentence ``$L_2(\cax/k,s)$ has a pole of order $t$ in $s=2$'' meaning
$$
\lim_{\Re(s)>2,s\to2}(s-2)^tL_2(\cax/k,s)=\alpha\ne0.
$$
Moreover, if $f(s)$ is a holomorphic function in $\Re(s)>\lambda$ and $\lim_{s\to\lambda}(s-\lambda)f(s)=\alpha\ne0$, we
will call $\alpha$ the residue of the function $f(s)$ in $s=\lambda$ and we will write $\res_{s=\lambda}(f(s))=\alpha$.
\end{enumerate}
\end{remark}

The Galois group $G_k:=\gal(\ov{k}/k)$ acts naturally on $\cag$. Let 
$\frO_{G_k}(\cag)$ be the set of $G_k$-orbits of $\cag$ with respect to this action.

Let $\Delta\subset\cac$ be the discriminant locus of $\phi$, i.e., the 
set of $y\in\cac$ such that the fiber $\phi^{-1}(y)=\caa_y$ is not smooth. Let $\cau:=\cac-\Delta$ be the smooth locus of $\phi$ and 
$\caa^{\sm}:=\phi^{-1}(\cau)$. 
Then $\phi$ induces a proper smooth morphism $\phi^{\sm}:\caa^{\sm}\to\cau$ whose fibers are abelian varieties. Let 
$\caf_a:=R^1(\phi^{\sm})_!\bbq_{\ell}$ and $\caf_b:=R^2(\phi^{\sm})_!\bbq_{\ell}$. These are 
smooth \'etale sheaves  on $\cau$ of pure weights 1 and 2, respectively. Let 
$\ov{\eta}$ be a geometric generic point of 
$\cau$ and $\pi_1(\cau,\ov{\eta})$ the algebraic fundamental group of $\cau$ with respect to 
$\ov{\eta}$. The stalks $(\caf_a)_{\ov{\eta}}$ and $(\caf_b)_{\ov{\eta}}$ of $\caf_a$ and $\caf_b$ at $\ov{\eta}$ 
are respectively isomorphic to $H^1(A)$ and $H^2(A)$. We denote by $\rho_a:\pi_1(\cau,\ov{\eta})\to\gl((\caf_a)
_{\ov{\eta}})$ and $\rho_b:\pi_1(\cau,\ov{\eta})\to\gl((\caf_b)_{\ov{\eta}})$ the monodromy representations 
associated to $\caf_a$, respectively $\caf_b$. Let $\caa':=\caa\times_{\cac}\cac'$.

\begin{theorem}\label{mainthm}
Suppose Tate's conjecture is true for $\caa'/k$ and that the monodromy representations $\rho_a$ and $\rho_b$ are 
irreducible. Then
\begin{equation}\label{rkthm}\begin{aligned}
\rk\left(\frac{A(K')}{\tau B(k)}\right)\le&\frac{\#\frO_{G_k}(\cag)}{\#\cag}(d\left(2d+1\right)(2g'-2))
\\&+\#\frO_{G_k}(\cag)2df_A.
\end{aligned}\end{equation}
\end{theorem}

\begin{remark}\label{remjac}
We will see in Remark \ref{remdim1} that when $A$ is a Jacobian variety then (\ref{rkthm}) can be improved to 
\begin{equation}\label{rkthm1}
\rk\left(\frac{A(K')}{\tau B(k)}\right)\le\frac{\#\frO_{G_k}(\cag)}{\#\cag}2d(2g'-2)+\#\frO_{G_k}(\cag)f_A,
\end{equation}
Moreover, the hypothesis that $\rho_b$ is irreducible is not necessary in this case.
\end{remark}

\begin{remark}\label{weil2}
If $A$ is an elliptic curve, then $\rho_a$ is automatically irreducible \cite[Lemme 3.5.5]{del81}.
\end{remark}

The following corollary is then a consequence of Theorem \ref{mainthm} and Remarks \ref{remjac} and \ref{weil2}.

\begin{corollary}\label{corec}
Suppose that $A$ is an elliptic curve and Tate's conjecture is true for $\caa'/k$. Then 
\begin{equation}\label{rkthma}
\rk(A(K'))\le\frac{\#\frO_{G_k}(\cag)}{\#\cag}(4g'-4)+\#\frO_{G_k}(\cag)f_A.
\end{equation}
\end{corollary}

\begin{remark}
We need the hypothesis that both $\rho_a$ and $\rho_b$ are irreducible in order to improve the trivial bound on the 
absolute values of the average traces of Frobenius (cf. Remark \ref{trivbd} and Section 3). There are proper flat morphisms 
$\psi:\cax\to\cac$
of relative dimension 1, where $\cax$ is a smooth projective irreducible surface, whose monodromy representations $\rho_a$
are irreducible (cf. \cite{kat81} and \cite[Th\'eor\`eme 2.3 and 2.4]{mic97}). So, if $\phi:\caa\to\cac$ is the associated Jacobian 
fibration, i.e., a 
proper flat morphism whose generic fiber is the Jacobian variety of the generic fiber of $\psi$, then the same holds for
the monodromy representation denoted also by $\rho_a$ associated to $\phi$. 

This improvement of the trivial bound as well as the use of the quite general 
Theorem \ref{thmhpw} are the main points in generalizing Silverman's result \cite[Theorem 13]{sil02} to the context of 
abelian varieties.
\end{remark}

\begin{remark}
One might naturally ask whether this theorem is true for covers which are geometrically Galois. 
For a more detailed discussion why we cannot 
expect to do so under the
method used in this paper to estimate the rank of the abelian variety see Remark \ref{galoisthm}.
\end{remark}

It will be shown in Proposition \ref{condprop} that $f_{A'}\le\#\cag f_A$ and that equality holds when $\pi$ is 
unramified. In this circumstance we rewrite (\ref{rkthm}), (\ref{rkthm1}) and (\ref{rkthma}) respectively as 
\begin{align}\label{rkthm3}
\rk\left(\frac{A(K')}{\tau B(k)}\right)&\le\frac{\#\frO_{G_k}(\cag)}{\#\cag}(d\left(2d+1)(2g'-2)
+2df_{A'}\right),
\\
\label{rkthm2}
\rk\left(\frac{A(K')}{\tau B(k)}\right)&\le\frac{\#\frO_{G_k}(\cag)}{\#\cag}(2d(2g'-2)+f_{A'})\text{ and }\\
\label{rkthmb}
\rk(A(K'))&\le\frac{\#\frO_{G_k}(\cag)}{\#\cag}(4g'-4+f_{A'}). 
\end{align}
Observe that (\ref{rkthmb}) is Silverman's result \cite[Theorem 13]{sil02}.

In the case where $\pi$ is unramified and $A$ is a Jacobian variety, as it happens in the case of elliptic curves, the 
new bound (\ref{rkthm2}) improves the geometric bound (\ref{ogg2}), if the action of $G_k$ on $\cag$ is non trivial. 
Observe that (\ref{rkthm}) does not improve (\ref{ogg2}) in general for an abelian variety
(even if $G_k$ acts non trivially on $\cag$).

Next we consider the problem of the variation of the rank of $A(K')/\tau B(k)$ along towers of function fields over $K$. More precisely,
we treat the following situation: 
$\cac$ is an elliptic curve without complex multiplication, 
$\cac'=\cac$ and $\pi=[n]:\cac\to\cac$ is the multiplication by $n$ map. In this case,
$\cag=\cac[n]$ is the subgroup of $n$-torsion points of $\cac$.
Denote $\cac_n:=\cac/[n]$ and $K_n:=k(\cac_n)$. 
The rank of $A(K_n)/\tau B(k)$ as $n\ge1$ varies grows slower than the geometric bound (cf. Theorem \ref{thmabtow}). 
This result is the exact analogue in the current situation of \cite[Theorem 16]{sil02}.
The reason why we consider this situation is that Serre's open image theorem 
\cite[Th\'eor\`eme $3'$]{serre72} 
implies that the action of $G_k$ on $\cag$ is highly non trivial as $n$ grows. See also Theorem \ref{jactow} for a 
discussion of the same problem in the case where $\cac$ has genus at least 2 and 
where $\cac'$ is the pull-back of $\cac$ under the multiplication by $n$ map $[n]:J_{\cac}\to J_{\cac}$ 
in the Jacobian variety $J_{\cac}$ of $\cac$. In this case, $\cag=J_{\cac}[n]$ (the $n$-torsion subgroup of $J_{\cac}$). 

In \S\,2 we describe the connection between Tate's conjecture and a conjecture related to the generalized analytic 
Nagao's conjecture. In \S\,3 we use local monodromy representations to give an upper bound for the absolute values of the average
traces of Frobenius in terms of the degree of the conductor of $A$. In \S\,4 we analyze how these upper bounds vary when we
extend the base $\cac$ by a finite geometrically abelian cover $\cac'\to\cac$ defined over $k$. In \S\,5 we prove Theorem \ref{mainthm}. 
In \S\,6 we give applications of Theorem \ref{mainthm}
to the rank of abelian varieties along towers of function fields over number fields.

\section{Tate's conjecture and a conjecture related to \\ the generalized Nagao's conjecture}

Given a prime ideal $\frp$ of $\cao_k$ and a smooth projective irreducible algebraic variety $\caz$ defined over $k$, we will denote by
${\caz}_{\frp}$ its reduction modulo $\frp$. 

Let $\psi:\cax\to\cay$ be a proper flat morphism of relative dimension $d$ of 
smooth projective connected algebraic varieties $\cax$ and $\cay$ defined as well as $\psi$ over $k$. Suppose that the 
generic fiber of $\psi$ is geometrically irreducible. Denote $m:=\dim(\cay)$. 

Let $S$ be a finite set of prime ideals of $\cao_k$ (which will be enlarged as needed). First we assume that for every
$\frp\notin S$, ${\cax}_{\frp}$ (resp. ${\cay}_{\frp}$) is a smooth projective connected variety over the
residue field $\bbf_{\frp}$ of $\frp$ of cardinality $q_{\frp}$.

For each $y\in{\cay}_{\frp}(\bbf_{\frp})$, let ${\cax}_{\frp,y}:={\psi}_{\frp}^{-1}(y)$ be the fiber of
${\psi}_{\frp}$ at $y$. Let ${F}_{\frp}$ be the topological generator of $\gal(\ov{\bbf}_{\frp}/\bbf_{\frp})$.
Denote also by $F_{\frp}$ its induced automorphism on $H^1({\cax}_{\frp,y})$.


For every $y\in{\cay}_{\frp}(\bbf_{\frp})$, let
$a_{\frp}({\cax}_{\frp,y}):=\tr({F}_{\frp}\,|\,H^1({\cax}_{\frp,y}))$ (respectively 
$b_{\frp}({\cax}_{\frp,y}):=\tr({F}_{\frp}\,|\,H^2({\cax}_{\frp,y}))$). The \emph{average traces of
Frobenius} are defined by

$$
\frA_{\frp}(\cax):=\frac1{q_{\frp}^m}\sum_{y\in\cay_{\frp}(\bbf_{\frp})}a_{\frp}({\cax}_{\frp,y}),
\text{ respectively }
\frB_{\frp}(\cax):=\frac1{q_{\frp}^m}\sum_{y\in{\cay}_{\frp}(\bbf_{\frp})}b_{\frp}({\cax}_{\frp,y}).
$$
Let $X/k(\cay)$ be the generic fiber of $\psi$, $P_X$ its Picard variety and $(\tau_X,Q_X)$ the $k(\cay)/k$-trace
of $P_X/k(\cay)$. 
Let $a_{\frp}(Q_X):=\tr(\Frob_{\frp}\,|\,H^1(Q_X)^{I_{\frp}})$. By base change (cf. \cite{mil80} or 
\cite[Appendix C]{har86})
 this number equals $\tr(F_{\frp}\,|\,H^1(Q_{X,\frp}))$. The \emph{reduced average trace of Frobenius} is 
defined 
by
$$
\frA_{\frp}^*(\cax):=\frA_{\frp}(\cax)-a_{\frp}(Q_X).
$$

\begin{theorem}\label{thmhpw}\cite[Th\'eor\`eme 1.2]{hinpacwa04}
Tate's Conjecture \ref{tateconj} for $\cax/k$ and $\cay/k$ implies Conjecture $M_{\text{an}}$:
\begin{multline*}
\res_{s=1}\left(\sum_{\frp\notin S}-\frA_{\frp}^*(\cax)\frac{\log(q_{\frp})}{q_{\frp}^s}\right)
+\res_{s=2}\left(\sum_{\frp\notin S}\frB_{\frp}(\cax)\frac{\log(q_{\frp})}{q_{\frp}^s}\right)\\
=\rk\left(\frac{P_X(k(\cay))}{\tau_XQ_X(k)}\right)+\rk(\ns(P_X/k(\cay))).
\end{multline*}
\end{theorem}

In the particular case where the relative dimension $d$ of $\psi$ is equal to 1, the Conjecture $M_{\text{an}}$ reduces
to the generalized analytic Nagao conjecture \cite[Remarque 1.4]{hinpacwa04}
\begin{equation}\label{gennagao}
\res_{s=1}\left(\sum_{\frp\notin S}-\frA_{\frp}^*(\cax)\frac{\log(q_{\frp})}{q_{\frp}^s}\right)=
\rk\left(\frac{P_X(k(\cay))}{\tau_XQ_X(k)}\right)
\end{equation}

In the hypothesis of Theorem \ref{mainthm}, Tate's conjecture is trivially true for the curve $\cac/k$. So, we need only
to assume its truth for $\caa'/k$.

\section{Estimates of average sums}

\begin{remark}\label{trivbd}
Let $\Delta_{\frp}$ be the discriminant locus of the reduction $\phi_{\frp}:\caa_{\frp}\to\cac_{\frp}$ of $\phi$ modulo $\frp$ for 
$\frp\notin S$. Let $\cau_{\frp}:=\cac_{\frp}-\Delta_{\frp}$ be the smooth locus of $\phi_{\frp}$. Note that after enlarging $S$, if
necessary, we may assume that $\Delta_{\frp}$ is obtained from $\Delta$ by reducing modulo $\frp$.

It follows from Weil's theorem that for every
$y\in{\cau}_{\frp}(\bbf_{\frp})$ all eigenvalues of $F_{\frp}$ acting on
$H^1({\caa}_{\frp,y})$ (respectively $H^2({\caa}_{\frp,y})$) have absolute value $\sqrt{q_{\frp}}$ 
(respectively $q_{\frp}$), thus
$a_{\frp}({\caa}_{\frp,y})=O(\sqrt{q_{\frp}})$ (respectively $b_{\frp}({\caa}_{\frp,y})=O(q_{\frp})$). For 
$y\in{\Delta}_{\frp}(\bbf_{\frp})$ it is a result due to
Deligne \cite[Th\'eor\`eme 3.3.1]{del81} that all the eigenvalues of $F_{\frp}$ acting on
$H^1({\caa}_{\frp,y})$ (respectively $H^2({\caa}_{\frp,y})$) have absolute value at most 
$\sqrt{q_{\frp}}$ (respectively $q_{\frp}$). So once again
$a_{\frp}({\cax}_{\frp,y})=O(\sqrt{q_{\frp}})$ (respectively $b_{\frp}({\cax}_{\frp,y})=O(q_{\frp})$). As a 
consequence, we have the trivial bound 
$\frA_{\frp}(\caa)=O(\sqrt{q_{\frp}})$ (respectively $\frB_{\frp}(\caa)=O(q_{\frp})$). Our goal in this 
paragraph is to use local monodromy representations to improve these estimates. The bounds obtained improve
even the ones coming from Deligne's equidistribution theorem (cf. \cite[Th\'eor\`eme 3.5.3]{del81} or 
\cite[(3.6.3)]{kat88}) which already improved the trivial bound (cf. Remark \ref{deleq}).
\end{remark}

The map ${\phi}_{\frp}$ restricts to a smooth proper morphism ${\phi}_{\frp}^{\sm}:{\caa}_{\frp}^{\sm}
\to{\cau}_{\frp}$ defined over $\bbf_{\frp}$. Let $\ell$ be a prime number different from the characteristic of $\bbf_{\frp}$. 
Let $\caf_{a,\frp}:=R^1({\phi}_{\frp}^{\sm})_!\bbq_{\ell}$ and $\caf_{b,\frp}:=R^2({\phi}_{\frp}^
{\sm})_!\bbq_{\ell}$ be smooth \'etale sheaves on ${\cau}_{\frp}$ of pure weights 1 and 2, respectively. 
Let $\ov{\eta}_{\frp}$ be the geometric 
generic point of ${\cau}_{\frp}$, ${A}_{\frp}$ the generic fiber of 
${\phi}_{\frp}$. As a consequence of proper base change 
\cite[VI, Corollary 2.7]{mil80}, $(\caf_a)_{\ov{\eta}}\cong H^1(A)\cong H^1({A}
_{\frp})\cong(\caf_{a,\frp})_{\ov{\eta}_{\frp}}$ (resp. 
$(\caf_b)_{\ov{\eta}}\cong H^2(A)\cong H^2({A}
_{\frp})\cong(\caf_{b,\frp})_{\ov{\eta}_{\frp}}$). The diagram
$$
\begin{CD}
\pi_1(\cau,\ov{\eta})@>{\rho_a}>>\gl((\caf_a)_{\ov{\eta}})\\
@VVV{\vdots}\\
\pi_1(\cau_{\frp},\ov{\eta}_{\frp})@>{\rho_{a,\frp}}>>\gl((\caf_{a,\frp})_{\ov{\eta}_{\frp}})
\end{CD}
$$
can be completed with a vertical downward arrow in the right hand side by using the former proper base change isomorphism. 
In this sense, we 
view $\rho_{a,\frp}$ as the reduction of $\rho_a$ modulo $\frp$. The same applies to $\rho_b$ and $\rho_{b,\frp}$. 

\begin{hypothesis}
We assume throughout this section that the representations $\rho_a$ and $\rho_b$ are irreducible.
\end{hypothesis}

Consequently, by construction, $\rho_{a,\frp}$ and $\rho_{b,\frp}$ are also irreducible. 

Let 
$$
\frA_{\frp}^{\sm}(\caa):=\frac1{q_{\frp}}\sum_{y\in{\cau}_{\frp}(\bbf_{\frp})}a_{\frp}({\caa}_{\frp,y})
\text{ and }
\frB_{\frp}^{\sm}(\caa):=\frac1{q_{\frp}}\sum_{y\in{\cau}_{\frp}(\bbf_{\frp})}b_{\frp}({\caa}_{\frp,y}).
$$
Observe that if $y\in{\cau}_{\frp}(\bbf_{\frp})$ then $y$ is identified with the map $y:\spec(\bbf_{\frp})\to{
\cau}_{\frp}$ which induces a homomorphism $\pi_1(\spec(\bbf_{\frp}))\cong G_{\bbf_{\frp}}\to\pi_1({\cau}_{\frp},\ov{\eta}_
{\frp})$ and the image by this map of $F_{\frp}$ is well defined up to conjugacy. The corresponding conjugacy class is
denoted by $F_{\frp,y}$. Note also that $y^*\caf_{a,\frp}\cong H^1({\caa}_{\frp,y})$ and $y^*\caf_{b,\frp}\cong H^2({\caa}_{\frp,y})$. 
Let $\tr(F_{\frp,y}\,|\,\caf_{a,\frp}):=\tr(\rho_{a,\frp}(F_{\frp,y}))$ and 
$\tr(F_{\frp,y}\,|\,\caf_{b,\frp}):=\tr(\rho_{b,\frp}(F_{\frp,y}))$. It is tautologous that 
$\tr(\rho_{a,\frp}(F_{\frp,y}))=\tr(F_{\frp}\,|\,y^*\caf_{a,\frp})=a_{\frp}({\caa}_{\frp,y})$ and 
$\tr(\rho_{b,\frp}(F_{\frp,y}))=\tr(F_{\frp}\,|\,y^*\caf_{b,\frp})=b_{\frp}({\caa}_{\frp,y})$. 
In particular,
$$
q_{\frp}\frA_{\frp}^{\sm}(\caa)=\sum_{y\in{\cau}_{\frp}(\bbf_{\frp})}\tr(F_{\frp,y}\,|\,\caf_{a,\frp})
\text{ and }
q_{\frp}\frB_{\frp}^{\sm}(\caa)=\sum_{y\in{\cau}_{\frp}(\bbf_{\frp})}\tr(F_{\frp,y}\,|\,\caf_{b,\frp}).
$$

Recall that $A/K$ is non-isotrivial, hence $\cau\varsubsetneq\cac$. Moreover, for every $\frp\notin S$ the discriminant locus $\Delta_
{\frp}=\cac_{\frp}-\cau_{\frp}$ of $\phi_{\frp}$ is equal to the reduction modulo $\frp$ of the discriminant locus $\Delta=\cac-\cau$ 
of $\phi$. Thus, 
${\cau}_{\frp}\varsubsetneq{\cac}_{\frp}$. As a consequence $H^0_c({\cau}_{\frp}\times_{\bbf_{\frp}}\ov{\bbf}
_{\frp},\caf_{a,\frp})=0$ and $H^0_c({\cau}_{\frp}\times_{\bbf_{\frp}}\ov{\bbf}
_{\frp},\caf_{b,\frp})=0$. Furthermore, since $\rho_{a,\frp}$ and $\rho_{b,\frp}$ are irreducible, it follows that $\pi_1(
{\cau}_{\frp},\ov{\eta}_{\frp})$ has no coinvariants, consenquently, $H^2_c({\cau}_{\frp}\times_{\bbf_{\frp}}
\ov{\bbf}_{\frp},\caf_{a,\frp})=0$ and $H^2_c({\cau}_{\frp}\times_{\bbf_{\frp}}\ov{\bbf}
_{\frp},\caf_{b,\frp})=0$. Now, by Grothendieck-Lefschetz formula \cite[VI, Theorem 13.4]{mil80}, 
\begin{equation}\label{grothle}\begin{aligned}
\sum_{y\in{\cau}_{\frp}(\bbf_{\frp})}\tr(F_{\frp,y}\,|\,\caf_{a,\frp})&=-\tr(F_{\frp}\,|\,H^1_c({\cau}_{\frp}
\times_{\bbf_{\frp}}\ov{\bbf}_{\frp},\caf_{a,\frp}))\text{ and}\\
\sum_{y\in{\cau}_{\frp}(\bbf_{\frp})}\tr(F_{\frp,y}\,|\,\caf_{b,\frp})&=-\tr(F_{\frp}\,|\,H^1_c({\cau}_{\frp}
\times_{\bbf_{\frp}}\ov{\bbf}_{\frp},\caf_{b,\frp})).
\end{aligned}\end{equation}

The Tate twists $\caf_{a,\frp}(1/2)$ and $\caf_{b,\frp}(1)$ of $\caf_{a,\frp}$ and 
$\caf_{b,\frp}$ are lisse \'etale sheaves of pure weight 0 on $\cau_{\frp}$. Hence, by \cite[Th\'eor\`eme 3.3.1]{del81}, the sheaves 
$H^1_c({\cau}_{\frp}
\times_{\bbf_{\frp}}\ov{\bbf}_{\frp},\caf_{a,\frp}(1/2))$ and $H^1_c({\cau}_{\frp}
\times_{\bbf_{\frp}}\ov{\bbf}_{\frp},\caf_{b,\frp}(1))$ are lisse \'etale sheaves of mixed weight at most 
1 over $\spec(\bbf_{\frp})$. Denote by $H^1_c({\cau}_{\frp}\times_{\bbf_{\frp}}\ov{\bbf}_{\frp},
\caf_{a,\frp}(1/2))^{\wt\le0}$, 
respectively $H^1_c({\cau}_{\frp}\times_{\bbf_{\frp}}\ov{\bbf}_{\frp},\caf_{b,\frp}(1))^{\wt\le0}$, the part of 
$H^1_c({\cau}_{\frp}\times_{\bbf_{\frp}}\ov{\bbf}_{\frp},\caf_{a,\frp}(1/2)))$, respectively 
$H^1_c({\cau}_{\frp}\times_{\bbf_{\frp}}\ov{\bbf}_{\frp},\caf_{b,\frp}(1)))$, of weight at most 0. 
By (\ref{grothle}) and Deligne's theorem \cite[Th\'eor\`eme 3.3.1]{del81},
\begin{multline}\label{grothle1}
\left|\sum_{y\in{\cau}_{\frp}(\bbf_{\frp})}\tr(F_{\frp,y}\,|\,\caf_{a,\frp}(1/2))\right|\\
\le O(1)+
\codim(H^1_c({\cau}_{\frp}\times_{\bbf_{\frp}}\ov{\bbf}_{\frp},\caf_{a,\frp}(1/2))^{\wt\le0})\sqrt{q_{\frp}},
\end{multline}
\begin{multline}\label{grothle2}
\left|\sum_{y\in{\cau}_{\frp}(\bbf_{\frp})}\tr(F_{\frp,y}\,|\,\caf_{b,\frp}(1))\right|\\
\le O(1)+
\codim(H^1_c({\cau}_{\frp}\times_{\bbf_{\frp}}\ov{\bbf}_{\frp},\caf_{b,\frp}(1))^{\wt\le0})\sqrt{q_{\frp}}.
\end{multline}
Since, for every $y\in{\cau}_{\frp}(\bbf_{\frp})$ we have $\tr(F_{\frp,y},\caf_{a,\frp}(1/2))=(1/\sqrt{q_{\frp}})
a_{\frp}({\caa}_{\frp,y})$ and $\tr(F_{\frp,y},\caf_{b,\frp}(1))=(1/q_{\frp})
b_{\frp}({\caa}_{\frp,y})$, the latter inequalities are equivalent to 
\begin{align}\label{apbd}
|\frA_{\frp}^{\sm}(\caa)|&\le O(1/\sqrt{q_{\frp}})+
\codim(H^1_c({\cau}_{\frp}\times_{\bbf_{\frp}}\ov{\bbf}_{\frp},\caf_{a,\frp}(1/2))^{\wt\le0})\text{ and }\\
\label{bpbd}|\frB_{\frp}^{\sm}(\caa)|&\le O(1)+
\codim(H^1_c({\cau}_{\frp}\times_{\bbf_{\frp}}\ov{\bbf}_{\frp},\caf_{b,\frp}(1))^{\wt\le0})\sqrt{q_{\frp}}.
\end{align}

Observe that 
$$
\frA_{\frp}(\caa)=\frA_{\frp}^{\sm}(\caa)+\frac1{q_{\frp}}\sum_{y\in{\Delta}_{\frp}(\bbf_{\frp})}a_{\frp}({
\caa}_{\frp,y})=\frA_{\frp}^{\sm}(\caa)+O(\sqrt{q_{\frp}}).
$$
Similarly, $\frB_{\frp}(\caa)=\frB_{\frp}^{\sm}(\caa)+
O(1)$. So, in order to estimate the absolute values of $\frA_{\frp}(\caa)$ and $\frB_{\frp}(\caa)$ it suffices to compute
the previous codimensions. These codimensions are expressed in terms of the multiplicities of the conductor of an 
abelian variety $X$ over a function field which we now recall.

Let $\caz$ be a smooth projective curve defined over a perfect
field $l$ of characteristic $p\ge0$, $\call=l(\caz)$ its function
field and $X$ an abelian variety defined over $\call$. Let $\ell\ne p$ be a prime number, $X[\ell]$ the subgroup of $\ell$-torsion
points of $X$, $T_{\ell}(X)$ the $\ell$-adic Tate module of $X$ and
$V_{\ell}(X)=T_{\ell}(X)\otimes_{\bbz_{\ell}}\bbq_{\ell}$. Let
$\call^s$ be a separable closure of $\call$. For every $z\in\caz$ 
denote by $I_z\subset\gal(\call^s/\call)$ the inertia subgroup corresponding to $z$ (which is well defined
up to conjugation).

\begin{definition}\label{defcond}
The multiplicity of the conductor $\frF_X$ of $X$ at $z\in\caz$ is equal to a sum of two numbers, the tame multiplicity $\epsilon_z$
and the wild multiplicity $\delta_z$ of $\frF_X$ in $z$. The first number is defined as
$\epsilon_z:=\text{codim}(V_{\ell}(X)^{I_z})$, where $V_{\ell}(X)^{I_z}$ denotes the subspace of $V_{\ell}(X)$ 
which is fixed by the action of $I_z$. Let 
$N(X)_z$ be the N\'eron model of $X$ over $\spec(\cao_{z})$ and denote by 
$\can(X)_z^0$ the connected component of the special fiber
$\can(X)_z$ of $N(X)_z$. Let $\kappa_z$ be the residue field of $z$ and $\ov{\kappa}_z$ its algebraic closure. 
Then $\can(X)_z^0\times_{\kappa_z}\ov{\kappa}_z$ is an extension of an abelian variety $Y_z$ by a linear algebraic group $\mathbb{L}_z$
both defined over $\ov{\kappa}_z$. Let $a_z:=\dim(Y_z)$.
The linear algebraic group $\mathbb{L}_z$ decomposes as $\mathbb{L}_z=\bbg_{m,\ov{\kappa}_z}^{t_z}
\times\bbg_{a,\ov{\kappa}_z}^{u_z}$, where $\bbg_{m,\ov{\kappa}_z}$ and $\bbg_{a,\ov{\kappa}_z}$ denote respectively the multiplicative 
and additive group scheme of $\ov{\kappa}_z$. The non-negative integers $a_z$, $t_z$ and $u_z$ 
are called the 
abelian, reductive and unipotent ranks of $\can(X)_z^0$, respectively. Then $\epsilon_z=2u_z+t_z$. Recall that 
$\dim(X)=a_z+t_z+u_z$.

The second number is defined as follows. Let $l_z/\kappa_z$ be a finite
Galois extension such that $G_{l_z}:=\gal(\ov{l}_z/l_z)$ acts trivially on $X[\ell]$. Let
$G_z:=\gal(l_z/\kappa_z)$, so $X[\ell]$ can be regarded as a $G_z$-module. Let $P_z$ be the projective
$\bbz_{\ell}[G_z]$-module whose character is the Swan character of $G_z$. Then
$\delta_z:=\dim_{\bbf_{\ell}}(\text{Hom}_{\bbz_{\ell}[G_z]}(P_z,X[\ell]))$
(cf. \cite[\S1]{ogg67}). This is a non-negative 
integer and is in fact independent from the choice of $l_z$.
\end{definition}

Let $\call(X[\ell])$ be the finite extension of $\call$ defined by adjoining the coordinates of the $\ell$-torsion points
of $X$. It follows from \cite[Remark (2), p. 480]{sertate} that $\delta_z=0$ for every $z\in\caz$ if and only if the 
extension $\call(X[\ell])/\call$ is tamely ramified.

\begin{remark}\label{tamecond}
After adjoining possibly finitely many prime ideals $\frp$ to $S$, we may assume that the conductor $\frF_{A_{\frp}}$ of $A_{\frp}$ is the
reduction of the conductor $\frF_A$ of $A$ modulo $\frp$. In particular, both conductors have the same degree. 
Moreover, since $\text{char}(K)=0$, there will be no contribution of the 
wild part of the conductor of $\frF_{A_{\frp}}$ for every $\frp\notin S$. 

In another direction, observe that ${\cac}_{\frp}$ has the same genus as $\cac$ for any $\frp\notin S$, after 
eventually extending $S$ even further.
\end{remark}

\begin{lemma}\label{codimlem}
\begin{align*}
\codim(H^1_c({\cau}_{\frp}\times_{\bbf_{\frp}}\ov{\bbf}_{\frp},\caf_{a,\frp}(1/2))^{\wt\le0})&=2d(2g-2)+\sum_{y
\in{\Delta}_{\frp}}\epsilon_y\text{ and }\\
\codim(H^1_c({\cau}_{\frp}\times_{\bbf_{\frp}}\ov{\bbf}_{\frp},\caf_{b,\frp}(1))^{\wt\le0})&\le d(2d-1)(2g-2)+
(2d-1)
\sum_{y\in{\Delta}_{\frp}}\epsilon_y.
\end{align*}
\end{lemma}

\begin{proof}
First note that by Poincar\'e's duality and the fact that $H^2_c({\cau}_{\frp}\times_{\bbf_{\frp}}
\ov{\bbf}_{\frp},\caf_{a,\frp})=0$ and $H^2_c({\cau}_{\frp}\times_{\bbf_{\frp}}\ov{\bbf}
_{\frp},\caf_{b,\frp})=0$ we obtain $H^0_{\et}({\cau}_{\frp}\times_{\bbf_{\frp}}
\ov{\bbf}_{\frp},\caf_{a,\frp})=0$ and $H^0_{\et}({\cau}_{\frp}\times_{\bbf_{\frp}}\ov{\bbf}
_{\frp},\caf_{b,\frp})=0$. The injection
${j}_{\frp}:{\cau}_{\frp}\hookrightarrow{\cac}_{\frp}$ \cite[2.0.7]{kat88} induces an exact sequence 
\begin{align*}
0&\to\bigoplus_{y\in{\Delta}_{\frp}}(\caf_{a,\frp}(1/2)_{\ov{\eta}_{\frp}})^{I_y}\to H^1_c({\cau}_{\frp}
\times_{\bbf_{\frp}}\ov{\bbf}_{\frp},\caf_{a,\frp}(1/2))\\
&\to H^1_{\et}({\cac}_{\frp}\times_{\bbf_{\frp}}\ov{\bbf}_
{\frp},(j_{\frp})_*\caf_{a,\frp}(1/2))\to 0.
\end{align*}
Since the last term is a sheaf of pure weight 1 over $\spec(\bbf_{\frp})$ \cite[Theorem 3.3.1]{del81}, we conclude that
$$
H^1_c({\cau}_{\frp}\times_{\bbf_{\frp}}\ov{\bbf}_{\frp},\caf_{a,\frp}(1/2))^{\wt\le0}\cong\bigoplus_{y\in{
\Delta}_{\frp}}(\caf_{a,\frp}(1/2)_{\ov{\eta}_{\frp}})^{I_y}.
$$
In a similar vein, 
$$
H^1_c({\cau}_{\frp}\times_{\bbf_{\frp}}\ov{\bbf}_{\frp},\caf_{b,\frp}(1))^{\wt\le0}\cong\bigoplus_{y\in{
\Delta}_{\frp}}(\caf_{b,\frp}(1)_{\ov{\eta}_{\frp}})^{I_y}.
$$

Since $(\caf_{a,\frp})_{\ov{\eta}_{\frp}}\cong H^1(A)$, then for every $y\in\Delta_{\frp}$ we have $\dim(\caf_{a,
\frp}(1/2)_{\ov{\eta}_{\frp}}^{I_y})=2d-\epsilon_y$. 
Moreover, 
$\dim(H^1_c({\cau}_{\frp}\otimes_{\bbf_{\frp}}
\ov{\bbf}_{\frp},\caf_{a,\frp}(1/2)))=-\chi_c({\cau}_{\frp}\times_{\bbf_{\frp}}\ov{\bbf}_{\frp},\caf_{a,\frp}(1/2))
=-\rk(\caf_{a,\frp}(1/2))\chi_c({\cau}_{\frp}\times_{\bbf_{\frp}}\ov{\bbf}_{\frp})=2d(2g-2+s_A)$, where $\chi_c$ denotes 
the Euler-Poincar\'e characteristic with respect to the cohomology with compact support, 
since for $\frp\notin S$ there is no contribution from the Swan conductor of $\caf_{a,\frp}$ (cf. Remark 
\ref{tamecond}). 
 Thus 
\begin{align*}
\codim(H^1_c({\cau}_{\frp}\times_{\bbf_{\frp}}\ov{\bbf}_{\frp},\caf_{a,\frp}(1/2))^{\wt\le0})&=2d(2g-2)+\sum_{y\in
{\Delta}_{\frp}}(2d)-\sum_{y\in{\Delta}_{\frp}}(2d-\epsilon_y)\\
&=2d(2g-2)+\sum_{y\in{\Delta}_{\frp}}\epsilon_y.
\end{align*}

Note that $(\caf_{b,\frp})_{\ov{\eta}_{\frp}}\cong H^2(A)\cong\bigwedge^2H^1(A)$. For every $y\in\Delta_{\frp}$ we have $\bigwedge^2
(H^1(A)^{I_y})$ \linebreak 
$\subset(\bigwedge^2H^1(A))^{I_y}$. Hence, $\dim(\caf_{b,\frp}(1)_{\ov{\eta}
_{\frp}}^{I_y})\ge(2d-\epsilon_y)(2d-\epsilon_y-1)/2$. Furthermore there is no contribution from the Swan conductor of $\caf_{b,\frp}$.
Thus, $\dim(H^1_c({\cau}_{\frp}\times_{\bbf_{\frp}}\ov{\bbf}_{\frp},\caf_{b,\frp}(1)))=-\chi_c({\cau}
_{{\frp}}\times
_{\bbf{\frp}}\ov{\bbf}_{\frp},\caf_{b,\frp}(1))=-\rk(\caf_{b,\frp}(1))\chi_c({\cau}_{\frp}\times
_{\bbf_{\frp}}\ov{\bbf}_{\frp})
=d(2d-1)(2g-2+s_A)$. Therefore,
\begin{align*}
&\codim(H^1_c({\cau}_{\frp}\times_{\bbf_{\frp}}\ov{\bbf}_{\frp},\caf_{b,\frp}(1))^{\wt\le0})\\
&\le d(2d-1)(2g-2)+\sum_{y
\in{\Delta}_{\frp}}{d(2d-1)}-\sum_{y\in{\Delta}_{\frp}}{(d-\epsilon_y/2)(2d-\epsilon_y-1)}\\
&\le d(2d-1)(2g-2)+(2d-1/2)\sum_{y\in{\Delta}_{\frp}}\epsilon_y-(1/2)\sum_{y\in\Delta_{\frp}}\epsilon_y^2\\
&\le d(2d-1)(2g-2)+(2d-1)\sum_{y\in{\Delta}_{\frp}}\epsilon_y.
\end{align*}
\end{proof}

\begin{theorem}\label{absbd}
\begin{align*}
|\frA_{\frp}(\caa)|&\le 2d(2g-2)+f_A+O(1/\sqrt{q_{\frp}})\\
|\frB_{\frp}(\caa)|&\le\left(d(2d-1)(2g-2)+(2d-1)f_A\right)\sqrt{q_{\frp}}+O(1).
\end{align*}
\end{theorem}

\begin{proof}
This follows from Lemma \ref{codimlem}, (\ref{apbd}), (\ref{bpbd}) and the fact that 
$$
\sum_{y\in{\Delta}_{\frp}}\epsilon_y=f_{{A}_{\frp}}=f_A
$$
(cf. Remark \ref{tamecond}).
\end{proof}

\begin{remark}\label{deleq}
Theorem \ref{absbd} in fact improves not only the trivial estimates on the absolute values of $\frA_{\frp}(\caa)$ and 
$\frB_{\frp}(\caa)$ of Remark \ref{trivbd}, but also upper bounds obtained by means of Deligne's equidistribution 
theorem (cf. \cite[Th\'eor\`eme
3.5.3]{del81} or \cite[(3.6.3)]{kat88}). Indeed, from the latter theorem
follows the bounds
\begin{align*}
\left|\sum_{y\in{\cau}_{\frp}(\bbf_{\frp})}\tr(F_{\frp,y}\,|\,\caf_a(1/2))\right|&\le\dim(H^1_c({\cau}_{\frp}
\times_{\bbf_{\frp}}\ov{\bbf}_{\frp},\caf_{a,\frp}(1/2)))\sqrt{q_{\frp}}\\
&\le2d(2g-2+s_A)\sqrt{q_{\frp}}\\
\left|\sum_{y\in{\cau}_{\frp}(\bbf_{\frp})}\tr(F_{\frp,y}\,|\,\caf_b(1))\right|&\le\dim(H^1_c({\cau}_{\frp}
\times_{\bbf_{\frp}}\ov{\bbf}_{\frp},\caf_{b,\frp}(1)))\sqrt{q_{\frp}}\\
&\le d(2d-1)(2g-2+s_A)\sqrt{q_{\frp}}.
\end{align*}
A fortiori,
\begin{align}\label{equia}
|\frA_{\frp}(\caa)|&\le2d(2g-2+f_A)+O(1/\sqrt{q_{\frp}})\\
\label{equib}
|\frB_{\frp}(\caa)|&\le d(2d-1)(2g-2+f_A)\sqrt{q_{\frp}}+O(1).
\end{align}
Hence both bounds are improved by Theorem \ref{absbd}.
\end{remark}

The geometric bound (\ref{ogg2}) involves the degree $f_{A'}$ of the conductor $\frF_{A'}$ of $A'$. The bound of Theorem
\ref{mainthm} is obtained in terms of the degree $f_A$ of the conductor $\frF_A$ of $A$. Although we will not need the 
following inequality in general, only in the unramified case, it holds in general.

\begin{proposition}\label{condprop}
$f_{A'}\le\#\cag f_A$. Equality holds if $\pi$ is unramified.
\end{proposition}

\begin{proof}
It suffices to prove that for each $y'\in\cac'$ such that $\pi(y')=y\in\cac$ we have the inequality $\epsilon_{y'}
\le\epsilon_y$. Indeed,
$$
f_{A'}=\sum_{y'}\deg(y')\epsilon_{y'}\le\sum_y\left(\sum_{y'|y}f(y'|y)\right)\epsilon_y\deg(y)\le\#\cag f_{A},
$$
where $f(y'|y)$ denotes the inertia degree of $y'$ over $y$.

Let $n\ge1$ be an integer invertible in the local ring $\cao_y$ of $y$. The base change morphism $N(A\times_KK_y)
\times_{\spec(\cao_y)}\spec(\cao_{y'})\to N(A'\times_{K'}K'_{y'})$ induces an injective (homo)mor\-phism at the level of 
torsion subgroups $N(A\times_KK_y)_{\tor}\times_
{\spec(\cao_y)}\spec(\cao_{y'})\to N(A'\times_{K'}K'_{y'})_{\tor}$, regarding them as products of group schemes of type $(p,\cdots,p)$. 
In particular, this map is also injective when 
restricted to $n$-torsion. Hence it induces, at the level of the connected 
components of the special fibers, an injection $\can(A\times_KK_y)^0_n\times
_{\spec(\kappa_y)}\spec(\ov{\kappa}_{y'})\cong(\bbz/n\bbz)^{2a_{y}+t_y}\hookrightarrow\can(A'\times_{K'}K'_{y'})^0_n\times
_{\spec(\kappa_{y'})}\spec(\ov{\kappa}_{y'})\cong(\bbz/n\bbz)^{2a_{y'}+t_{y'}}$, i.e., $2a_y+t_y\le2a_{y'}+t_{y'}$. This 
latter inequality is equivalent to $\epsilon_{y'}\le\epsilon_y$.

If $\pi$ is unramified, then the base change morphism is already an isomorphism 
(cf. \cite[Chapter 7, Theorem 1, p. 176]{ray90}), in particular it also induces an 
isomorphism at the level of the connected component of the special fiber, a fortiori $\epsilon_{y'}=\epsilon_y$.
\end{proof}

\section{Estimates after base change}

\begin{lemma}\label{trfrob}
Let $\frp\notin S$ and ${\pi}_{\frp}:{\cac}'_{\frp}\to{\cac}_{\frp}$ be the reduction of $\pi$ modulo
$\frp$, $y'\in{\cac}'_{\frp}(\bbf_{\frp})$ and $y={\pi}_{\frp}(y')\in{\cac}_{\frp}(\bbf_{\frp})$. Then
$a_{\frp}({\caa}'_{\frp,y'})=a_{\frp}({\caa}_{\frp,y})$ and $b_{\frp}({\caa}'_{\frp,y'})=b_{\frp}(
{\caa}_{\frp,y})$.
\end{lemma}

\begin{proof}
This is immediate from the fact that ${\caa'}_{\frp,y'}\cong{\caa}_{\frp,y}$, once $y$ and $y'$ are defined
over $\bbf_{\frp}$.
\end{proof}

After adding finitely many prime ideals to $S$, we may also assume that for every $\frp\notin S$,
${\pi}_{\frp}:{\cac}'_{\frp}\to{\cac}_{\frp}$ is also finite geometric abelian with 
with geometric automorphism group 
${\cag}_{\frp}:=\aut_{\ov{\bbf}_{\frp}}({\cac}'_{\frp}/{\cac}_{\frp})$ isomorphic to $\cag$. Let
${\bbg}_{\frp}:=\aut_{\bbf_{\frp}}({\cac}'_{\frp}/{\cac}_{\frp})$.

Let $\cah$ be a subgroup of $\cag$, $\cac'':=\cac'/\cah$ the intermediate curve $\cac'\to\cac''\to\cac$ 
with $\aut_{\ov{k}}(\cac'/\cac'')=\cah$ and $\pi':\cac'\to\cac''$ the first covering. 
Since $\cag$ is abelian, the cover $\pi'':\cac''\to\cac$ is also geometrically Galois. 
Since $\cag$ is finite, there exists a finite extension $L$ of $k$ such that all elements of $\cag$ are defined over $L$. So 
$G_k$ acts on $\cag$ through the finite quotient $\gal(L/k)$.
Hence all the intermediate curves $\cac''$ are defined over $L$. If
$\gal(L/k)$ acts on $\cah$, then $\cac''$ is also defined over $k$.

Let $D$ be the reduced ramification divisor of $\pi:\cac'\to\cac$ and $D_{\frp}$ the reduced ramification divisor of $\pi_{\frp}:\cac'_
{\frp}\to\cac_{\frp}$. 
After extending $S$, if necessary, we may assume that $D_{\frp}$ is the reduction of $D$ modulo $\frp$. In particular, 
\begin{equation}\label{suppeq}
\#\supp({D}_{\frp})=\#\supp(D)\le\deg(D). 
\end{equation}
Let $\car_{\frp}:=\cau_{\frp}\cap\supp(D_{\frp})$ and $\cav_{\frp}:=\cau_{\frp}-\car_{\frp}$. 

The next proposition is proved exactly as in \cite[Proposition 11]{sil02}.

\begin{proposition}\label{countpt}
Let $F_{\frp}$ be a generator of $\gal(\ov{\bbf}_{\frp}/\bbf_{\frp})$ and $\cah\subset\cag$ the subgroup defined by
$\cah:=\{F_{\frp}(x)\circ x^{-1}\,;\,x\in\cag\}$.
\begin{enumerate}
\item The group $\cah$ is defined over $\bbf_{\frp}$, hence the curve ${\cac}''_{\frp}:={\cac}'_{\frp}/\cah$ is
also defined over $\bbf_{\frp}$.

Denote $\cau'_{\frp}:=\pi^{-1}_{\frp}(\cau_{\frp})$,
$\car'_{\frp}:=\pi^{-1}_{\frp}(\car_{\frp})$ and $\cav'_{\frp}:=\pi^{-1}_{\frp}(\cav_{\frp})=\cau'_{\frp}-\car'_{\frp}$. Denote similarly
$\cau''_{\frp}$, $\car''_{\frp}$ and $\cav''_{\frp}$ with respect to the map $\pi''_{\frp}$.

\item The image of ${\cav}'_{\frp}(\bbf_{\frp})$ in ${\cac}_{\frp}(\bbf_{\frp})$ coincides with the image of
${\cav}''_{\frp}(\bbf_{\frp})$ in ${\cac}_{\frp}(\bbf_{\frp})$.

\item For every $y\in{\cav}_{\frp}(\bbf_{\frp})$ such that $y={\pi}_{\frp}(y')$ with $y'\in
{\cav}'_{\frp}(\bbf_{\frp})$ there exist exactly $\#{\bbg}_{\frp}$ points in 
${\cav}'_{\frp}(\bbf_{\frp})$ lying over $y$.

\item For every $y\in{\cac}_{\frp}(\bbf_{\frp})$ such that $y={\pi}''_{\frp}(y'')$ with
$y''\in{\cav}''_{\frp}(\bbf_{\frp})$ there exist exactly $(\cag:\cah)$ points of
${\cav}''_{\frp}(\bbf_{\frp})$ lying over $y$.
\end{enumerate}
\end{proposition}

\begin{proposition}(cf. \cite[Proposition 12]{sil02})\label{avgtr2}
\begin{align*}
|\frA_{\frp}(\caa')|&\le\frac{\#{\bbg}_{\frp}}{\#\cag}(2d(2g'-2))+\#{\bbg}_{\frp}f_{A}+O(1/\sqrt{q_{\frp}})\\
|\frB_{\frp}(\caa')|&\le\sqrt{q_{\frp}}\left(\frac{\#{\bbg}_{\frp}}{\#\cag}d(2d-1)(2g'-2)+\#{\bbg}_{\frp}
(2d-1)f_A\right)+O(1).
\end{align*}
\end{proposition}

\begin{proof}
Let 
$$
\frA_{\frp}^{\sm,\ur}(\caa):=\frac1{q_{\frp}}\sum_{y\in{\cav}_{\frp}(\bbf_{\frp})}a_{\frp}({\caa}_{\frp,y}).
$$
So, by Weil's theorem and (\ref{suppeq}),
\begin{align*}
|\frA_{\frp}^{\sm}(\caa)-\frA_{\frp}^{\sm,ur}(\caa)|&=\frac1{q_{\frp}}\left|\sum_{y\in{\car}_{\frp}(\bbf_{\frp})}
a_{\frp}({\caa}_{\frp,y})\right|\le\frac1{q_{\frp}}2d\sqrt{q_{\frp}}\#\supp({D}_{\frp})\\
&\le2d(1/\sqrt{q_{\frp}})\deg(D)=O(1/\sqrt{q_{\frp}}).
\end{align*}
Since $\frA_{\frp}(\caa)=\frA_{\frp}^{\sm}(\caa)+O(1/\sqrt{q_{\frp}})$, we conclude that $\frA_{\frp}(\caa)=\frA_{\frp}
^{\sm,\ur}(\caa)+O(1/\sqrt{q_{\frp}})$.

Let 
$$
\frA_{\frp}^{\sm}(\caa'):=\frac1{q_{\frp}}\sum_{y'\in{\cau}_{\frp}'(\bbf_{\frp})}a_{\frp}({\caa}_{\frp,y'}')
\text{ and }
\frA_{\frp}^{\sm,\ur}(\caa'):=\frac1{q_{\frp}}\sum_{y'\in{\cav}_{\frp}'(\bbf_{\frp})}a_{\frp}({\caa}_
{\frp,y'}').
$$
Note that $\#\car_{\frp}'\le\deg(\pi)\#\car_{\frp}\le\deg(\pi)\deg(D)$. 
As above,
\begin{align*}
|\frA_{\frp}^{\sm}(\caa')-\frA_{\frp}^{\sm,\ur}(\caa')|&=\frac1{q_{\frp}}\left|\sum_{y'\in{\car}'_{\frp}(\bbf_
{\frp})}a_{\frp}({\caa}'_{\frp,y'})\right|\\
&\le2d(1/\sqrt{q_{\frp}})\deg(D)\deg(\pi)=O(1/\sqrt{q_{\frp}}).
\end{align*}
Since $\frA_{\frp}(\caa')=\frA_{\frp}^{\sm}(\caa')+O(1/\sqrt{q_{\frp}})$, then $\frA_{\frp}(\caa')=\frA_{\frp}^{\sm,\ur}
(\caa')+O(1/\sqrt{q_{\frp}})$. Define similarly $\frA_{\frp}^{\sm}(\caa'')$ and $\frA_{\frp}^{\sm,\ur}(\caa'')$.

As a consequence of Lemma \ref{trfrob} and Proposition \ref{countpt} we obtain 
\begin{align*}
q_{\frp}\frA_{\frp}^{\sm,\ur}(\caa')&=\sum_{y'\in{\cav}'_{\frp}(\bbf_{\frp})}a_{\frp}({\caa}'_{\frp,y'})\\
&=\#{\bbg}_{\frp}\sum_{y\in{\pi}_{\frp}({\cav}_{\frp}'(\bbf_{\frp}))}
a_{\frp}({\caa}_{\frp,y})=\#{\bbg}_{\frp}\sum_{y\in{\pi}_{\frp}''({\cav}_{\frp}''(\bbf_{\frp}))}
a_{\frp}({\caa}_{\frp,y})\\
&=\frac{\#{\bbg}_{\frp}}{(\cag:\cah)}\sum_{y''\in{\cav}_{\frp}''(\bbf_{\frp})}a_{\frp}({\caa}
_{\frp,y''}'')=\frac{\#{\bbg}_{\frp}}{(\cag:\cah)}q_{\frp}\frA_{\frp}^{\sm,\ur}(\caa'').
\end{align*}
Thus, by Theorem \ref{absbd}, 
\begin{align*}
|\frA_{\frp}^{\sm,\ur}(\caa')|&=\frac{\#{\bbg}_{\frp}}{(\cag:\cah)}|\frA_{\frp}^{\sm,\ur}(\caa'')|\le
\frac{\#{\bbg}_{\frp}}{(\cag:\cah)}(2d(2g''-2)+f_{A''})+O(1/\sqrt{q_{\frp}})
\\ &\le\frac{\#{\bbg}_{\frp}}{\#\cag}2d(2g'-2)+\#{\bbg}_{\frp}f_A+O(1/\sqrt{q_{\frp}}),
\end{align*}
where the latter inequality follows from the Riemann-Hurwitz formula applied to $\pi':\cac'\to\cac''$ and and Proposition 
\ref{condprop}. Note that if $\pi$ is unramified, then the latter inequality is actually an equality. 
Whence the same 
inequality holds for $|\frA_{\frp}(\caa')|$.

Let 
$$
\frB_{\frp}^{\sm,\ur}(\caa):=\frac1{q_{\frp}}\sum_{y\in{\cav}_{\frp}(\bbf_{\frp})}b_{\frp}({\caa}_{\frp,y}).
$$
Similarly, $\frB_{\frp}^{\sm,\ur}(\caa)=\frB_{\frp}(\caa)+O(1)$ and $\frB_{\frp}^{\sm,\ur}(\caa')=
({\#{\bbg}_{\frp}}/{(\cag:\cah)})
\frB_{\frp}^{\sm,\ur}(\caa'')$ and therefore (again by Theorem \ref{absbd})
\begin{align*}
|\frB_{\frp}^{\sm,\ur}(\caa')|&=\frac{\#{\bbg}_{\frp}}{(\cag:\cah)}|\frB_{\frp}^{\sm,\ur}(\caa'')|\\
&\le\frac{\#{\bbg}_{\frp}}{(\cag:\cah)}\sqrt{q_{\frp}}
(d(2d-1)(2g''-2)+(2d-1)f_{A''})
+O(1)\\
&\le \sqrt{q_{\frp}}\left(\frac{\#{\bbg}_{\frp}}{\#\cag}d(2d-1)(2g'-2)+\#{\bbg}_{\frp}
(2d-1)f_A\right)+O(1),
\end{align*}
where again in the latter inequality we used the Riemann-Hurwitz formula and Proposition \ref{condprop}. 
Whence the same inequality holds for $|\frB_
{\frp}(\caa')|$.
\end{proof}

\begin{remark}[Towards Galois covers]\label{galoisthm}
The reason why we need to introduce the intermediate curve $C'/{\cah}$ in order to obtain an estimate for $|\frA_{\frp}(\caa')|$, and so
have to suppose that the cover $\cac'\to\cac$ is geometrically abelian instead of geometrically Galois, is that there may be points in 
$\cac'_{\frp}$ lying over an
$\bbf_{\frp}$-rational point of $\cac_{\frp}$ which are not $\bbf_{\frp}$-rational. This is related to curves which are twistings of 
$\cac'_{\frp}$ over $\bbf_{\frp}$. We now introduce those curves.

For each $\sigma\in{\bbg}_{\frp}$ there exists a smooth projective irreducible curve $\cac'_{\frp,\sigma}$ defined over $\bbf_{\frp}$ 
whose function field $\bbf_{\frp}(\cac'_{\frp,\sigma})$ equals the subfield of $\bbf_{\frp^m}(\cac'_{\frp})$ which is fixed by 
$\sigma F_{\frp}$, where $\bbf_{\frp^m}$ denotes the extension of $\bbf_{\frp}$ of degree $m$.
Clearly, $\cac_{\frp,\text{id}}=\cac'_{\frp}$.
We have the well-known Twisting Lemma \cite{bom}
$$
\#{\cac}_{\frp}(\bbf_{\frp})=\frac1{\#{\bbg}_{\frp}}\sum_{\sigma\in{\bbg}_{\frp}}\#{\cac}_{\frp,\sigma}'(\bbf_
{\frp}).
$$
More precisely, for each $y\in{\cac}_{\frp}(\bbf_{\frp})$ we have $\#{\bbg}_{\frp}$ points in the curves ${\cac}_{\frp,\sigma}'$ 
that are rational over $\bbf_{\frp}$ and lie over $y$. 
Indeed, given $y'\in\cac'_{\frp}$ such that $y'\mapsto y$, then there exists $\sigma\in\bbg_{\frp}$
such that $F_{\frp}y'=\sigma^{-1}y'$, i.e., $y'\in\cac'_{\frp,\sigma}(\bbf_{\frp})$. Note that $\sigma'\in\bbg_{\frp}$ satisfies 
$(\sigma')^{-1}y'=\sigma^{-1}y'$ if and only if $\sigma(\sigma')^{-1}$ 
lies in the decomposition group of $y'$. Let $n_y$ be the number of 
elements in the orbit of $y'$ with respect to $\bbg_{\frp}$. So, the number of $\bbf_{\frp}$-rational points of the curves $\cac'
_{\frp,\sigma}$ which lie above $y$ is equal to $n_y$ times $\#\bbg_{\frp}/n_y$, i.e., $\#\bbg_{\frp}$.

Each twisting $\cac_{\frp,\sigma}'$ comes equipped with a Galois cover 
$\pi_{\frp,\sigma}:\cac_{\frp,\sigma}'\to\cac_{\frp}$ defined over $\bbf_{\frp}$ 
with group $\bbg_{\frp}$. Let $\caa'_{\frp,\sigma}:=\caa_{\frp}\times_{\cac_{\frp}}\cac_{\frp,\sigma}'$. Let
$$
\frA_{\frp,\sigma}(\caa'):=\frac1{q_{\frp}}\sum_{y_{\sigma}'\in\cac_{\frp,\sigma}'(\bbf_{\frp})}a_{\frp}(\caa_{\frp,\sigma,y_{\sigma}'}').
$$
We apply Lemma \ref{trfrob} to the twistings $\cac_{\frp,\sigma}'$ of $\cac'_{\frp}$ 
and obtain (from the previous discussion on the twisting lemma)
\begin{equation}\label{twistequ}
\sum_{\sigma\in\bbg_{\frp}}\frA_{\frp,\sigma}(\caa')=\#\bbg_{\frp}\frA_{\frp}(\caa).
\end{equation}

Theorem \ref{thmhpw} implies an upper bound of the rank of $A(K)/\tau B(k)$ in terms of the local data given by the average traces of 
Frobenii $\frA_{\frp}^*(\caa)$ and $\frB_{\frp}(\caa)$. The characteristic $p$ information $\frA_{\frp,\sigma}(\caa')$, if $\sigma\ne
\text{id}$, does not come in the expression of an upper bound for the rank of $A(K_{\sigma}')/\tau B(k)$ in the same way, 
where $K_{\sigma}'=k(\cac'_
{\sigma})$ and $\cac'_{\sigma}$ is a twisting of $\cac'$ defined over $k$. 
One reason for this is that as $\frp$ varies so does $\sigma$. Hence, even under the identification of $\cag$ 
and $\cag_{\frp}$ for $\frp\notin S$, there is no ``general'' $\sigma\in\cag$ that would give a twist $\cac'_{\sigma}$ of $\cac'$ defined
over $k$ that would reduce modulo all prime ideals $\frp\notin S$ to $\cac'_{\frp,\sigma}$. If that were true then the inequality
$$
\rk\left(\frac{A(K')}{\tau B(k)}\right)\le\sum_{\sigma}\rk\left(\frac{A(K_{\sigma}')}{\tau B(k)}\right)
$$
together with (\ref{twistequ}) would give a proof of our result for geometrically Galois covers $\cac'\to\cac$ (cf. \S5).
\end{remark}

\section{Proof of Theorem \ref{mainthm}}

\begin{proof}
It follows from Theorem \ref{thmhpw} that 
\begin{align*}
&\rk\left(\frac{A(K')}{\tau B(k)}\right)\le\rk\left(\frac{A(K')}{\tau B(k)}\right)+\rk(\ns(A/K))\\
&=\res_{s=1}\left(\sum_{\frp\notin S}-\frA_{\frp}^*(\caa')\frac{\log(q_{\frp})}{q_{\frp}^s}\right)+
\res_{s=2}\left(\sum_{\frp\notin S}\frB_{\frp}(\caa')\frac{\log(q_{\frp})}{q_{\frp}^s}\right).
\end{align*}
Note that
$$
|\frA_{\frp}^*(\caa')|\le|\frA_{\frp}(\caa')|+|a_{\frp}(B)|\le|\frA_{\frp}(\caa')|+2\dim(B)\sqrt{q_{\frp}}.
$$
So, by Proposition \ref{avgtr2},
\begin{equation}\label{rkbd1}\begin{aligned}
&\rk\left(\frac{A(K')}{\tau B(k)}\right)
\le\left(\frac{2d(2g'-2)}{\#\cag}+f_A\right)\res_{s=1}\left(\sum_{\frp\notin S}{\#{\bbg}_{\frp}}
\frac{\log(q_{\frp})}{q_{\frp}^s}\right)\\
&+O\left(\res_{s=1}\left(\sum_{\frp\notin S}\frac{\log(q_{\frp})}{q_{\frp}^{s+1/2}}\right)\right)+2\dim(B)\res_{s=1}
\left(\sum_{\frp\notin S}\frac{\log(q_{\frp})}{q_{\frp}^{s-1/2}}\right)\\
&+\left(\frac{d(2d-1)(2g'-2)}{\#\cag}+(2d-1)f_A\right)\res_{s=2}\left(\sum_{\frp\notin S}{\#{\bbg}_{\frp}}
\frac{\log(q_{\frp})}{q_{\frp}^{s-1/2}}\right)\\
&
+O\left(\res_{s=2}\left(\sum_{\frp\notin S}\frac{\log(q_{\frp})}{q_{\frp}^{s}}\right)\right).
\end{aligned}\end{equation}
The third series converges for $\Re(s)>1/2$, thus the corresponding residue equals 0. The fourth series converges for
$\Re(s)>3/2$ and can be extended meromorphically to the whole plane with just a simple pole at $s=3/2$. 
Hence, the fourth residue also equals 0. The last residue is also equal to 0,
since the series converges for $\Re(s)>1$. 

As we have already observed, we have fixed a finite Galois extension $L$ of $k$ such that $G_k$ acts on $\cag$ via $\gal(L/k)$. 
Moreover, for each $\sigma$ in the $\frp$-th Frobenius conjugacy class $(\frp,L/k)$, the group $\bbg_{\frp}$ corresponds bijectively 
to the elements of $\{x\in\cag\,|\,\sigma(x)=x\}$. We denote the cardinality of this set by $h^0(\sigma,\cag)$. By (\ref{rkbd1})
\begin{equation}\label{rkbd3}\begin{aligned}
&\rk\left(\frac{A(K')}{\tau B(k)}\right)\\
&\le\left(\frac{2d(2g'-2)}{\#\cag}+f_A\right)\left(\sum_{\sigma\in\gal(L/k)}{h^0(\sigma,\cag)}
\res_{s=1}\left(\sideset{}{'}\sum_{\frp}\frac{\log(q_{\frp})}{q_{\frp}^s}\right)\right)\\
&+\left(\frac{d(2d-1)(2g'-2)}{\#\cag}+(2d-1)f_A\right)\times\\
&\left(\sum_{\sigma\in\gal(L/k)}{h^0(\sigma,\cag)}
\res_{s=2}\left(\sideset{}{'}\sum_{\frp}\frac{\log(q_{\frp})}{q_{\frp}^{s-1/2}}\right)\right),
\end{aligned}\end{equation}
where $\sum_{\frp}'$ denotes the sum over all $\frp\notin S$ such that $\sigma\in(\frp,L/k)$. 
It follows from \cite[Corollary of Theorem 8, \S4, Chapter VIII]{lan70} (cf. also 
the proof of \cite[Proposition 4]{sil00}) that all two residues are equal to $1/\#\gal(L/k)$. The result
now follows from \cite[Lemma 9]{sil02}.

The second statement follows from the next remark.
\end{proof}

\begin{remark}\label{remdim1}
In the case $A$ is the Jacobian variety of the generic fiber of a proper flat morphism $\phi:\cax\to\cac$ of relative 
dimension 1, where $\cax$ is a smooth projective irreducible surface defined over $k$, 
then Conjecture M$_{\text{an}}$ (Theorem \ref{thmhpw}) reduces to the generalized Nagao conjecture 
(\ref{gennagao}) \cite[Remarque 3.4]{hinpacwa04}. Hence we do not need to consider the estimates on $|\frB_{\frp}(\caa)|$ and the 
hypothesis of the irreducibility of $\rho_b$ is no longer necessary.
So, (\ref{rkbd3}) implies (\ref{rkthm1}), hence (\ref{rkthma}).
\end{remark}

\section{The rank in special towers}

In this paragraph we apply Theorem \ref{mainthm} to the special tower described in the Introduction. 

\begin{theorem}[Serre, \cite{serre72}, Th\'eor\`eme $3'$]\label{serrethm}
Let $\cac/k$ be an elliptic curve defined over a number field $k$ without complex multiplication. 
There is an integer $I(\cac,k)$ depending only on $\cac$ and $k$ so that for every 
integer $n\ge1$ the image of the Galois representation $\rho_{\cac,n}:G_k\to\aut(\cac[n])\cong\gl_2(\bbz/n\bbz)$ has index at 
most $I(\cac,k)$ in $\aut(\cac[n])$.
\end{theorem}

The same proof as in \cite[Theorem 16]{sil02} yields the following result.

\begin{theorem}\label{thmabtow}
Let $\cac$ be an elliptic curve without complex multiplication defined over a number field $k$. Suppose that $\cac'=\cac$ and 
$\pi$ is the multiplication
by $n$ map $[n]:\cac\to\cac$. Let $\cac_n:=\cac/[n]$, $K_n:=k(\cac_n)$ its function field and $\caa_n:=\caa\times_{\cac}\cac_n$.
Assume that Tate's conjecture is true for $\caa_n/k$ and the monodromy
representations $\rho_a$ and $\rho_b$ are irreducible. Then
\begin{enumerate}
\item For every integer $n\ge1$:
$$
\rk\left(\frac{A(K_n)}{\tau B(k)}\right)\le I(\cac/k)\frac{d(n)}{n^2}2df_{A_n},
$$
if $d\ge2$, where $d(n)$ denotes the number of positive divisors of $n$. If $d=1$, we may replace $2d$ by $1$.
\item The sum
$$
\frac1x\sum_{n\le x}\frac1{\log(f_{A_n})}\rk\left(\frac{A(K_n)}{\tau B(k)}\right)
$$
is bounded as $x\to\infty$. Thus the average rank of $A(K_n)/\tau B(k)$ is smaller than a fixed multiple of the
logarithmic of the degree $f_A$ of the conductor $\frF_A$ of $A$.
\item There exists a constant $\kappa=\kappa(k,\cac,A)$ 
so that for sufficiently large $n$ we have
$$
\rk\left(\frac{A(K_n)}{\tau B(k)}\right)\le f_{A_n}^{\kappa/\log(\log(f_{A_n}))}.
$$
In particular, for every $\epsilon>0$ we have
$$
\rk\left(\frac{A(K_n)}{\tau B(k)}\right)\ll f_{A_n}^{\epsilon},
$$
where the implied constant depends on $k$, $\cac$, $A$ and $\epsilon$, but not on $n$.
\end{enumerate}
\end{theorem}

We now consider the situation where $\cac$ has genus at least 2 and 
$\cac'$ is the pull-back of $\cac$ under the multiplication by $n$ map 
$[n]:J_{\cac}\to J_{\cac}$ in the 
Jacobian variety of $\cac$ and $\pi$ is the corresponding unramified geometrically abelian cover $\pi:\cac'\to\cac$. In this case, $\cag=
J_{\cac}[n]$. In order to state our
result we make a short digression on the Mumford-Tate group and the Mumford-Tate conjecture in the case of abelian varieties.

Let $X$ be an abelian variety defined over a number field $k$ of dimension $d\ge1$. Let $\ell$ be a prime number, $X[\ell]$ the subgroup 
of $\ell$-torsion points of $X$, $T_{\ell}(X)$ the $\ell$-adic Tate group of $X$ and $V_{\ell}(X):=T_{\ell}(X)\times_{\bbz_{\ell}}\bbq_
{\ell}$. Let $X^{\vee}:=\pic^0(X)$ be the dual abelian variety of $X$. After eventually extending $k$ and making an isogeny, we may 
assume that we have a degree 1 polarization $\lambda:X\to X^{\vee}$. Let $n\ge1$ be an integer, $\mu_{\ell^n}$ the group of $\ell^n$-th
roots of unity in $\mathbb{C}$ and $\mu_{\ell^{\infty}}:=\bigcup_{n\ge1}\mu_{\ell^n}$. The
polarization $\lambda$ yields non-degenerate alternating (Weil) pairings $X[\ell]\times X[\ell]\to\mu_{\ell}$ and $T_{\ell}(X)\times T_
{\ell}(X)\to\mu_{\ell^{\infty}}$. The natural action of $G_k$ on $X[\ell]$, hence on $T_{\ell}(X)$, and $\lambda$ induce 
Galois representations $\rho_{\ell}:G_k\to\Aut(X[\ell])\cong\gsp_{2d}(\bbf_{\ell})$ and $\rho_{\ell^{\infty}}:G_k\to\Aut(V_{\ell}(X))\cong
\gsp_{2d}(\bbq_{\ell})$.

Let $G_{k,\ell^{\infty}}:=\rho_{\ell^{\infty}}(G_k)$ and $G_{k,\ell}:=\rho_{\ell}(G_k)$. Let 
$$
\rho_{\infty}:=\prod_{\ell}\rho_{\ell^{\infty}}:G_k\to\Aut(V_{\ell}(X))\cong\prod_{\ell}\gsp_{2d}(\bbq_{\ell}).
$$
Serre proves in 
\cite[Th\'eor\`eme $1'$]{serre86cf}, \cite{serre86rib}, \cite{serre86mf} 
 that if $k$ is large enough, depending on $X$, then $\rho_{\infty}(G_k)$ is an open subgroup of 
$\prod_{\ell}G_{k,\ell^{\infty}}$. 

Fix an isomorphism $V_{\ell}(X)\cong\bbq_{\ell}^{2d}$. Let $\und{\gl}_{2d,\bbq_{\ell}}$ be the linear algebraic group defined in 
dimension $2d$ over $\bbq_{\ell}$. Denote by $\und{G}_{k,\ell^{\infty}}$ the Zariski closure of $G_{k,\ell^{\infty}}$ on $\und{\gl}
_{2d,\bbq_{\ell}}$ and $\und{G}_{k,\ell^{\infty}}^0$ its identity component.

The Mumford-Tate group $\mt(X)$ of $X$ is defined as follows. Fix an embedding $k\hookrightarrow\bbc$. 
Let $V:=H^1(X,\bbq)$ be the first singular 
cohomology group of $X$ with rational coefficients. A complex structure on $V_{\bbr}:=H^1(X,
\bbr)$ is an $\bbr$-linear map $J:H^1(X,\bbr)\to H^1(X,\bbr)$ such that $J^2=-I$, where $I$ denotes the identity. Let $S^1=\{z\in\bbc\,|\,
|z|=1\}$. To give $J$ is equivalent to giving a representation of real algebraic groups $h:S^1\to\gl(V_{\bbr})$ defined by $a+ib\mapsto
aI+bJ$. Fix an isomorphism $V\cong\bbq^{2d}$. 
The group $\mt(X)$ is defined as the smallest algebraic subgroup $G$ of $\und{\gl}_{2d,\bbq}$ defined over $\bbq$ such that 
$h(S^1)\subset G(\bbr)$. 

The Mumford-Tate conjecture states that $\und{G}_{k,\ell^{\infty}}^0=\mt(X)_{\ell}:=\mt(X)\times_{\bbq}\bbq_{\ell}$. R. Pink proved in 
\cite[Theorem 5.4]{pin98} under the hypothesis that $\End_{\ov{k}}(X)=\bbz$ and numerical conditions on $d=\dim(X)$ that both groups are
equal to $\und{\gsp}_{2d,\bbq_{\ell}}$. Since $A$ has a degree 1 polarization $\lambda$, 
$\und{G}_{k,\ell^{\infty}}$ is an algebraic subgroup
of $\und{\gsp}_{2d,\bbq_{\ell}}$. So in the latter circumstance $\und{G}_{k,\ell^{\infty}}=\und{G}_{k,\ell^{\infty}}^0=\und{\gsp}_{2d,
\bbq_{\ell}}$. Serre had previously proved this theorem under the 
assumptions that $\End_{\ov{k}}(X)=\bbz$, $d$ is odd or 2 or 6 \cite[Th\'eor\`eme 3]{serre86cf}. For other cases where the 
Mumford-Tate is known see \cite{pin98}.

Bogomolov's algebraicity theorem states that $G_{k,\ell^{\infty}}$ has finite index in the subgroup $\und{G}_{k,\ell^{\infty}}
(\bbq_{\ell})$ of $\bbq_{\ell}$-rational points of $\und{G}_{k,\ell^{\infty}}$ 
\cite[2.2.1]{serre85cf}.
In particular, if $\und{G}_{k,\ell^{\infty}}=\und{G}_{k,\ell^{\infty}}^0=\und{\gsp}_{2d,\bbq_{\ell}}$, then 
\begin{enumerate}
\item[(a)] $G_{k,\ell^{\infty}}$ is an open subgroup of the subgroup $\gsp_{2d}(\bbq_{\ell})$ of $\bbq_{\ell}$-rational points of 
$\und{\gsp}_{2d,\bbq_{\ell}}$. 
\end{enumerate}
Furthermore, it is also known that the rank $r$ of $\und{G}_{k,\ell^{\infty}}$ does not depend on $\ell$ \cite[2.2.4]{serre85cf} 
and that if $r=d+1$ and $\End_{\ov{k}}(X)=\bbz$, 
then $\und{G}_{k,\ell^{\infty}}=\und{G}_{k,\ell^{\infty}}^0=\und{\gsp}_{2d,\bbq_{\ell}}$ \cite[2.2.7]{serre85cf}
and 
\begin{enumerate}
\item[(b)] for almost all $\ell$, $G_{k,\ell}$ coincides with $\gsp_{2d}(\bbf_{\ell})$ \cite[Th\'eor\`eme 3]{serre86mf} which is equal to 
the reduction $\mt(X)(\ell)$ of $\mt(X)_{\ell}$ modulo $\ell$. 
\end{enumerate}
So in the same way as in \cite[Introduction]{serre72} (a) and (b) imply the following theorem.

\begin{theorem}\label{avopim}
Let $X/k$ be an abelian variety defined over a number field. Suppose that $\End_{\ov{k}}(X)=\bbz$
and $r=d+1$, then there exists an integer $I(X,k)$ depending only on $X$ and $k$ such that the 
image of the Galois representation $\rho_{X,n}:G_k\to\aut(X[n])\cong\gsp_{2d}(\bbz/n\bbz)$ has index at most $I(X,k)$ in $\aut(X[n])$.
\end{theorem}

\begin{remark}
As observed by Serre the previous condition on $r$ is satisfied if $d$ is odd or 2 or 6 \cite[2.2.8]{serre85cf}.
\end{remark}

\begin{theorem}\label{jactow}
Let $\cac_n:=\cac'$ be the pull-back of $\cac$ under the multiplication by $n$ map 
$[n]:J_{\cac}\to J_{\cac}$ in the 
Jacobian variety of $\cac$ and $\pi$ the corresponding unramified geometrically abelian cover $\pi:\cac'\to\cac$. Let $K_n:=k(\cac_n)$ be 
the function field of $\cac_n$, $\caa_n:=\caa\times_{\cac}\cac_n$, $A_n/K_n$ the generic fiber of $\caa_n\to\cac_n$, $g$, resp. $g_n$, 
the genus of $\cac$, resp. $\cac_n$. Let $\und{G}_{k,\ell^{\infty}}$ be the Zariski closure of the image of $\rho_{\ell^{\infty}}:G_k
\to\aut(V_{\ell}(J_{\cac}))\cong\gsp_{2g}(\bbq_{\ell})$. 
Assume that Tate's conjecture is true for $\caa_n/k$, the monodromy
representations $\rho_a$ and $\rho_b$ are irreducible, $\End_{\ov{k}}(J_{\cac})=\bbz$ and the rank $r$ of $\und{G}_{k,\ell^{\infty}}$ is 
equal to $g+1$. Then
\begin{enumerate}
\item For every $n\ge1$ we have
\begin{equation}\label{est3}
\rk\left(\frac{A(k_n)}{\tau B(k)}\right)\le
I(J_{\cac},k)\frac{d(n)}{n^{2g}}(d(2d+1)(2g_n-2)+2df_{A_n}).
\end{equation}
If $d=1$, we may replace $d(2d+1)$ by $2d$ in the first sumand of the latter expression and $2d$ by 1 in the second summand.
\item The sum
$$
\frac1{x}\sum_{n\le x}\frac1{\log(f_{A_n})}\rk\left(\frac{A(K_n)}{\tau B(k)}\right)
$$
is bounded as $x\to\infty$. Thus the average rank of $A(K_n)/\tau
B(k)$ is smaller than a fixed multiple of the logarithm of the degree
$f_A$ of the conductor $\frF_A$ of $A$.
\item There exists a constant $\kappa=\kappa(k,J_{\cac},A)$ 
so that for sufficiently large $n$ we have
$$
\rk\left(\frac{A(K_n)}{\tau B(k)}\right)\le f_{A_n}^{\kappa/\log\log(f_{A_n})}.
$$
In particular, for every $\epsilon>0$ we have
$$
\rk\left(\frac{A(K_n)}{\tau B(k)}\right)\ll f_{A_n}^{\epsilon},
$$
where the implied constant depends only on $k$, $J_{\cac}$, $A$ and
$\epsilon$, but not on $n$.
\end{enumerate}
\end{theorem}

\begin{proof}
By Theorem \ref{avopim} and \cite[Lemma 10]{sil02} 
\begin{equation}\label{orbit1}\begin{aligned}
\#\frO_{G_k}(J_{\cac}[n])&\le
I(J_{\cac},k)\#\frO_{\aut(J_{\cac}[n])}(J_{\cac}[n])\\
&=I(J_{\cac},k)\#\frO_{\gsp_{2g}(\bbz/n\bbz)}
((\bbz/n\bbz)^{2g}).
\end{aligned}\end{equation}
The same proof as in \cite[Proposition 15]{sil02} replacing the
statement that every non zero vector in a finite dimensional vector space can be extended
to a basis of the vector space by Witt's theorem \cite[Chapter XIV,
Theorem 2, p. 360]{lan65} shows that the number of orbits of the
action of $\gsp_{2g}(\bbz/n\bbz)$ in $(\bbz/n\bbz)^{2g}$ is also equal
to $d(n)$. Hence, by (\ref{orbit1}), 
\begin{equation}\label{orbit2}
\#\frO_{G_k}(J_{\cac}[n])\le I(J_{\cac},k)d(n).
\end{equation}
Item (1) now follows from Theorem \ref{mainthm} and
(\ref{orbit2}). 

Notice that it
follows from \cite[Theorem 3.3]{apo76} that 
$$
\sum_{n\le x}d(n)\sim x\log(x).
$$
Therefore,
\begin{equation}\label{avg1}
\frac1x\sum_{2\le n\le x}\frac{d(n)}{\log(n)}
\end{equation}
is bounded for every $x\ge2$. 
Moreover, by Proposition \ref{condprop}, $f_{A_n}=n^{2g}f_A$, and by
the Riemann-Hurwitz formula, since $\pi$ is unramified,
$2g_n-2=n^{2g}(2g-2)$. Thus, by (1),
$$
\rk\left(\frac{A(K_n)}{\tau B(k)}\right)\le I(J_{\cac},k)d(n)(d(2d+1)(2g-2)+2df_A)
$$ 
and  item (2) follows from (\ref{avg1}).

By \cite[Theorem 13.12]{apo76}
\begin{equation}\label{ap2}
\limsup_{n\to\infty}\frac{\log(d(n))}{\log(n)/\log\log(n)}=\log(2). 
\end{equation}
In particular, by (1), there exist absolute constants
$c_1,c_2,c_3,c_4$ such that 
\begin{align*}
\rk\left(\frac{A(K_n)}{\tau B(k)}\right)&\le
c_1I(J_{\cac},k)(d(2d+1)(2g-2)+2df_A)n^{c_2/\log\log(n)}\\
&\le c_3I(J_{\cac},k)(d(2d+1)(2g-2)+2df_A)f_{A_n}^{c_4/\log\log(f_{A_n})},
\end{align*}
which proves item (3).
\end{proof}

\end{document}